\documentclass[12pt,leqno]{amsart}
\usepackage{amssymb}
\usepackage{amsmath,amssymb,color}
\begin{document}
\baselineskip=18pt

%
\newtheorem{thm}{Theorem}[section]
\newtheorem{lem}[thm]{Lemma}
\newtheorem{cor}[thm]{Corollary}
\newtheorem{prop}[thm]{Proposition}
\newtheorem{defn}[thm]{Definiction}
\newtheorem{rmk}{Remark}
\newtheorem{exa}{Example}
\numberwithin{equation}{section}
\newcommand{\DDDD}{D_{\rho,\delta}}
\newcommand{\weight}{e^{2s\varphi}}
\newcommand{\tilweight}{e^{2s\widetilde{\va}}}
\newcommand{\ep}{\varepsilon}
\newcommand{\la}{\lambda}
\newcommand{\va}{\varphi}
\newcommand{\ppp}{\partial}
\newcommand{\izt}{\int^t_0}
\newcommand{\ddd}{\mathcal{D}}
\newcommand{\www}{\widetilde}
\newcommand{\pppg}{\partial_{t,\gamma}^{\alpha}}
\newcommand{\pppa}{\partial_t^{\alpha}}
\newcommand{\VVVV}{Cs^3e^{Cs}\mathcal{D}(u)^2 + Cs^3e^{2s\sigma_1}M^2}
\newcommand{\ddda}{d_t^{\alpha}}
\newcommand{\DDDa}{D_t^{\alpha}}
\newcommand{\lla}{L^{\frac{3}{4}}}
\newcommand{\llb}{L^{-\frac{1}{4}}}
\newcommand{\HHMIN}{H^{-1}(\Omega)}
\newcommand{\HHONE}{H^1_0(\Omega)}
\newcommand{\sumij}{\sum_{i,j=1}^n}
\newcommand{\Halp}{H_{\alpha}(0,T)}
\newcommand{\Hal}{H^{\alpha}(0,T)}
\newcommand{\lllb}{L^{-\frac{3}{4}}}
\newcommand{\pdif}[2]{\frac{\partial #1}{\partial #2}}
\newcommand{\ppdif}[2]{\frac{\partial^2 #1}{{\partial #2}^2}}
\newcommand{\uen}{u_N^{\ep}}
\newcommand{\fen}{F_N^{\ep}}
\newcommand{\penk}{p^{\ep}_{N,k}}
\newcommand{\penl}{p^{\ep}_{N,\ell}}
\newcommand{\R}{\mathbb{R}}
\newcommand{\Q}{\mathbb{Q}}
\newcommand{\Z}{\mathbb{Z}}
\newcommand{\C}{\mathbb{C}}
\newcommand{\N}{\mathbb{N}}
\newcommand{\uu}{\mathbf{u}}
\renewcommand{\v}{\mathbf{v}}
\newcommand{\y}{\mathbf{y}}
\newcommand{\RR}{\mathbf{R}}
\newcommand{\Y}{\mathbf{Y}}
\newcommand{\w}{\mathbf{w}}
\newcommand{\z}{\mathbf{z}}
\newcommand{\G}{\mathbf{G}}
\newcommand{\ooo}{\overline}
\newcommand{\OOO}{\Omega}
\newcommand{\CCO}{{_{0}C^1[0,T]}}
\newcommand{\dddaO}{{_{0}\ddda}}
\newcommand{\WWW}{{_{0}W^{1,1}(0,T)}}
\newcommand{\CCOO}{{^{0}C^1[0,T]}}
\newcommand{\HH}{_{0}{H^{\alpha}}(0,T)}
\newcommand{\HHO}{_{0}H^1(0,T)}
\newcommand{\hhalf}{\frac{1}{2}}
\newcommand{\DDD}{\mathcal{D}}
\newcommand{\RRR}{\mathcal{R}}
\newcommand{\RRRR}{\longrightarrow}
\newcommand{\lbra}{_{H^{-1}(\Omega)}{\langle}}
\newcommand{\wdel}{\www{\delta}}
\newcommand{\weps}{\www{\varepsilon}}
\newcommand{\rbra}{\rangle_{H^1_0(\Omega)}}
\newcommand{\lbran}{_{(H^{-1}(\Omega))^N}{\langle}}
\newcommand{\rbran}{\rangle_{(H^1_0(\Omega))^N}}
\newcommand{\llll}{L^{\infty}(\Omega\times (0,t_1))}
\newcommand{\sumk}{\sum_{k=1}^{\infty}}
\newcommand{\OOOPM}{\Omega_{\pm}}
\renewcommand{\baselinestretch}{1.5}
\renewcommand{\div}{\mathrm{div}\,}  
\newcommand{\grad}{\mathrm{grad}\,}  
\newcommand{\rot}{\mathrm{rot}\,}  
\newcommand{\xxx}{x', x_n}
\newcommand{\xxz}{x', 0}
\newcommand{\urur}{\left( \frac{1}{r}, \theta\right)}

\allowdisplaybreaks
\renewcommand{\theenumi}{\arabic{enumi}}
\renewcommand{\labelenumi}{(\theenumi)}
\renewcommand{\theenumii}{\alph{enumii}}
\renewcommand{\labelenumii}{(\theenumii)}
\def\thefootnote{{}}

\title
[]
{Inverse parabolic problems of determining functions
with one spatial-component
independence by Carleman estimate
}


\author{
$^1$ O.~Yu.~Imanuvilov,
$^2$ Y. Kian and \, $^{3,4,5}$ M.~Yamamoto }

\thanks{
$^1$ Department of Mathematics, Colorado State
University, 101 Weber Building, Fort Collins, CO 80523-1874, U.S.A.
e-mail: {\tt oleg@math.colostate.edu}
\\
$^2$ Aix Marseille Univ, Universit\'e de Toulon, CNRS, CPT,
Marseille, France
e-mail:{\tt yavar.kian@univ-amu.fr}
\\
$^3$ Graduate School of Mathematical Sciences, The University
of Tokyo, Komaba, Meguro, Tokyo 153-8914, Japan
\\
$^4$ Honorary Member of Academy of Romanian Scientists,
Splaiul Independentei Street, no 54,
050094 Bucharest Romania
\\
$^5$ Peoples' Friendship University of Russia
(RUDN University) 6 Miklukho-Maklaya St, Moscow, 117198, Russian Federation
e-mail: {\tt myama@ms.u-tokyo.ac.jp}
}

\date{}
\maketitle

\baselineskip 18pt

\begin{abstract}
For an initial-boundary value problem for a parabolic
equation in the spatial variable $x=(x_1,.., x_n)$ and time
$t$, we consider an inverse problem of determining
a coefficient which is independent of one spatial component $x_n$
by extra lateral boundary data.
We apply a Carleman estimate to prove a conditional stability
estimate for the inverse problem.  Also we prove similar results
for the corresponding inverse source problem.
\\
{\bf Key words.}
inverse source problem, inverse coefficient problem,
Carleman estimates, stability
\\
{\bf AMS subject classifications.}
35R30, 35R25
\end{abstract}

\section{Introduction and the main results}
Let $x = (x',x_n) \in \R^n$ with $x' = (x_1,..., x_{n-1})\in \R^{n-1}$ be
the spatial variable, $t$ the time variable,
and $D \subset \R^{n-1}$ be a bounded domain with smooth boundary $\ppp D$,
and $\ell>0$ a constant.
We set 
$$
\OOO = D \times (0,\ell).
$$
We note that $\OOO \subset \R^n$ is a cylindrical domain of spatial
variables $x \in \R^n$.

We consider
$$
\ppp_tv = \Delta v + p(x',t)v(x,t), \quad 
(x,t)\in \OOO\times (-\delta,\delta)          \eqno{(1.1)}
$$
with $p \in L^{\infty}(D \times (-\delta, \delta))$.
For notational convenience, we choose $-\delta < t < \delta$ as
time interval with initial time $-\delta$, not $0$.
Henceforth we denote $\ppp_{x_i} = \frac{\ppp}{\ppp x_i}$,
$\ppp_{x_i}\ppp_{x_j} = \frac{\ppp^2}{\ppp x_i\ppp x_j}$, $\ppp_{x_i}^2
= \frac{\ppp^2}{\ppp x_i^2}$, $1\le i,j \le n$,
$\ppp_t = \frac{\ppp}{\ppp t}$, $\nabla = (\ppp_{x_1}, ..., \ppp_{x_n})$,
$\nabla_{x,t} = (\nabla, \ppp_t)$, 
$\Delta = \sum_{i=1}^n \ppp_{x_i}^2$.
We set
$$
H^{2,1}(\OOO \times (-\delta, \delta))
= \{ v=v(x,t);\, v, \ppp_tv, \ppp_{x_i}v, \ppp_{x_i}\ppp_{x_j}v 
\in L^2(\OOO\times (-\delta, \delta)), \, 1\le i,j\le n\}.
$$

Our method allows us to consider more general parabolic equations and even
hyperbolic equations.
However for demonstrating the key idea, it is sufficient to consider (1.1).
\\

In equation (1.1), the coefficient $p=p(x',t)$ is assumed to be independent of
the one component $x_n$ of the spatial variables, 
but dependent on time, which describes
that a physical property under consideration
may change in time and independent of the depth variable $x_n$
for example.

The main purpose of this article is to establish the stability
for the inverse coefficient problem of determining
$p=p(x',t)$ by extra data $\ppp_{x_n} v$ on
a lateral subboundary $\Gamma \times (0,\ell) \times (-\delta, \delta)$,
where $\Gamma \subset \ppp D$ is a given subboundary.

For studying the inverse coefficient problem, we first consider
an inverse source problem for a parabolic
equation:
$$\left\{ \begin{array}{rl}
&\ppp_tu(x',x_n,t) = \Delta u + p_0(x',t)u + R(x',x_n,t)f(x',t),
\quad (x,t)\in \OOO\times (-\delta,\delta), \\
& u(\xxz,t) = \ppp_{x_n}u(\xxz,t) = 0, \quad (x',t)\in D\times (-\delta,\delta),
\end{array}\right.
                                    \eqno{(1.2)}
$$
where $p_0 \in L^{\infty}(D \times (-\delta,\delta))$ is given.

In (1.2) we assume
$$
R \in C^1([0,\ell]; L^{\infty}(D \times (-\delta,\delta)))
                                              \eqno{(1.3)}
$$
and
$$
R(x',0,t) \ne 0, \quad x' \in \ooo{D}, \, -\delta \le t \le \delta.
                                        \eqno{(1.4)}
$$

Let $\Gamma \subset \ppp D$ be an arbitrarily fixed non-empty relatively
open subset.  We arbitrarily choose a subdomain $D_0 \subset D$
satisfying
$$
\left\{ \begin{array}{rl}
& \ooo{D_0} \subset D \cup \Gamma, \quad
\mbox{$\ppp D_0 \cap \ppp D$ is a non-empty relatively open subset
of $\ppp D$}, \\
& \mbox{and $\ooo{\ppp D_0 \cap \ppp D} \subset \Gamma$.}
\end{array}\right.                    \eqno{(1.5)}
$$

Here and henceforth $\ooo{D}$ denotes the closure of a set $D$.

The inverse source problem is formulated as follows:
{\it determine a factor $f(x',t)$ of the source term by $u, \nabla u$ 
on $\Gamma \times (0,\ell) \times (-\delta,\delta)$.}
\\

Moreover we assume that
$u, \ppp_{x_n}u \in H^{2,1}(\OOO \times (-\delta, \delta))$ and 
$$
\Vert \ppp_{x_n}u(\cdot,\delta)\Vert_{H^1(\OOO)}
+ \Vert \ppp_{x_n}u(\cdot,-\delta)\Vert_{H^1(\OOO)}
+ \sum_{k=0}^1 \Vert \nabla_{x,t}^k \ppp_{x_n}u\Vert_{L^2(\ppp D \times
(0, \ell) \times (-\delta,\delta))}
$$
$$
+ \Vert \ppp_{x_n}u\Vert^2_{L^2(-\delta,\delta;H^2(\ppp\OOO))}
+ \sum_{k=0}^1 \Vert \nabla_{x,t}\ppp_{x_n}u(\cdot,\ell,\cdot)
\Vert_{L^2(D\times (-\delta,\delta))} \le M,                     \eqno{(1.6)}
$$
where a positive constant $M>0$ is arbitrarily fixed a priori bound.
We set
$$
\mathcal{D}(u) = \left( \int_{\Gamma\times (0,\ell) \times (-\delta,\delta)}
(\vert \nabla_{x,t} \ppp_{x_n}u\vert^2
+ \vert \ppp_{x_n}u\vert^2) d\sigma dt
+ \Vert \ppp_{x_n} u \Vert^2_{L^2(-\delta,\delta;
H^2(\Gamma\times (0,\ell))} \right)^{\frac{1}{2}}.         \eqno{(1.7)}
$$
We have
\\
{\bf Theorem 1.1.}\\
{\it Let $u$ satisfy equation 
(1.2), $p_0 \in L^{\infty}(D \times (-\delta,\delta))$ and conditions 
(1.3), (1.4) and (1.6) hold true. 
\\
(i) For a given subdomain $D_0 \subset D$ satisfying (1.5), 
there exist some constants
$0 < \delta_0 < \delta$, $\theta \in (0,1)$ and $C>0$ such that
$$
\Vert f\Vert_{L^2(D_0\times (-\delta_0,\delta_0))}
\le C\mathcal{D}(u)^{\theta}.
$$
\\
(ii) For any given $\delta_1 \in (0, \delta)$, there exist constants
$C_1 > 0$, $\theta_1 \in (0,1)$ and a subdomain $D_1 \subset D$ satisfying
(1.5) such that
$$
\Vert f\Vert_{L^2(D_1\times (-\delta_1,\delta_1))}
\le C_1\mathcal{D}(u)^{\theta_1}.
$$}
\\

In terms of the function $d \in C^2(\ooo{D})$ constructed in
Lemma 2.1 in Section 2, we can rewrite Theorem 1.1 (i) as follows:
If 
$$
\delta_0 < \left( \frac{\min_{x'\in \ooo{D_0}} d(x')}
{\max_{x'\in \ooo{D}} d(x')} \right)^{\frac{1}{2}}\delta,   \eqno{(1.8)}
$$
then the conclusion of (i) holds.

We note that
if $D_0$ is larger, that is, if we want to determine $f$ in a larger
spatial subdomain $D_0$, then
$\min_{x'\in \ooo{D_0}} d(x')$ is smaller, so that
the time interval $(-\delta_0,\delta_0)$ when we can prove a stability 
estimate  for the function $f$, becomes
smaller.  Moreover, since $d\vert_{\ppp\OOO\setminus \Gamma}
= 0$ by (2.1) stated in Section 2, 
the quantitiy $\min_{x'\in \ooo{D_0}} d(x')$ in (1.8) 
tends to $0$ if $D_0$ approaches to $D$, that is, we cannot estimate 
$f$ in $D$ even if we choose any short time interval. 
In Theorem 1.1 (ii), as is seen by the proof below, when $\delta_1 < \delta$ is
closer to $\delta$, the stability subdomain $D_1$ becomes smaller.

To sum up, Theorem 1.1 asserts a conditional stability estimate in determining
a factor $f$ of the source term in a proper subdomain of the cylinder $D 
\times (-\delta, \delta)$, which
holds conditionally with an a priori boundedness condition (1.6).

Theorem 1.1 (i) yields a uniqueness result only for $x' \in D$ at $t=0$.
More precisely we state
\\
{\bf Corollary 1.1.}\\{\it Let $u$ satisfy equation
(1.2), $p_0 \in L^{\infty}(D \times (-\delta,\delta))$ and conditions
(1.3), (1.4) and (1.6) hold true. 
If
$$
\ppp_{x_n}u = \nabla_{x'}\ppp_{x_n}u = 0 \quad \mbox{on 
$\Gamma \times (0,\ell) \times
(-\delta, \delta)$},
$$
then
$$
f(x',0) = 0, \qquad x' \in D.
$$}
\\

So far, we do not know the uniqueness in the whole spatial domain 
$D$ even on very small time interval.
\\

Next we state the main result on the inverse coefficient problem.
For $p = p(x',t)$ and $g_0,g_1\in H^{\frac{3}{2},\frac{3}{4}}
(D\times(-\delta,\delta))$, let $v = v(p)(x,t)$ satisfy
$$\left\{\begin{aligned}
\ppp_tv = \Delta v + p(x',t)v(x,t),& \quad (x,t) \in \OOO\times 
(-\delta,\delta),   \\
\partial_{x_n}^kv(x',0,t)=g_k(x',t),&\quad x'\in D, \, -\delta < t < \delta,
\quad k=0,1.
\end{aligned}\right.
                      \eqno{(1.9)}
$$

Fixing $g_0$ and $g_1$ in $D \times (-\delta,\delta)$, we estimate 
a coefficient $p(x',t)$ by data norm of the solutions 
defined by (1.7).  We can change $g_0$ and $g_1$, but we fix them,  
which simplifies and does not affect the essence of the arguments.   
\\

For the statement of our main result, we introduce an admissible set of
unknown coefficients $p(x',t)$.
For arbitrarily fixed constant $M>0$, we define an admissible set
$\mathcal{P}$ by
\begin{align*}
& \mathcal{P} = \{ p(x',t);\,
\Vert p\Vert_{L^{\infty}(D\times (-\delta,\delta))} \le M, 
 v(p), \, \ppp_{x_n}v(p) \in H^{2,1}(\OOO \times (-\delta,\delta))
\cap L^{\infty}(\OOO \times (-\delta,\delta)), \\
& \mbox{$v(p)$ satisfies (1.9)}\}.
\end{align*}

We are ready to state a conditional stability estimate for the inverse
coefficient problem.
\\
{\bf Theorem 1.2.}\\{\it
(i) For a given subdomain $D_0 \subset D$ satisfying (1.5),
there exist constants
$\delta_0 \in (0,\delta)$, $\theta \in (0,1), \alpha_0>0$ and $C>0$ such that
$$
\Vert p-q\Vert_{L^2(D_0\times (-\delta_0,\delta_0))}
\le C\mathcal{D}(v(p) - v(q))^{\theta}
$$
for $p, q \in \mathcal{P}$ if
$$
\mbox{either} \quad 
\vert v(p)(x',0,t)\vert \ge \alpha_0 \quad \mbox{for all 
$(x',t)\in \ooo{D_0}\times [-\delta,\delta]$} 
$$
$$ 
\quad \mbox{or} \quad \vert v(q)(x',0,t)\vert \ge \alpha_0 
\quad \mbox{for all $(x',t)\in \ooo{D_0}\times [-\delta,\delta]$}. 
                                                \eqno{(1.10)}
$$
\\
(ii) For any given $\delta_1 \in (0, \delta)$, there exist constants
$C_1 > 0$, $\theta_1 \in (0,1)$ and a subdomain $D_1 \subset D$ satisfying
(1.5) such that
$$
\Vert p-q\Vert_{L^2(D_1\times (-\delta_1,\delta_1))}
\le C_1\mathcal{D}(v(p) - v(q))^{\theta_1}
$$
for $p, q \in \mathcal{P}$ if (1.10) holds for $D_0=D_1.$}
\\

An inverse problem of determining an $x_n$-independent
factor of the source term, is considered in Beznoshchenko \cite{Bez1, Bez2},
Gaitan and Kian \cite{GK}, Kian and Yamamoto \cite{KY}.
Our proof is different from \cite{Bez1}, \cite{Bez2}, \cite{KY}, and
based on a Carleman estimate.
The work \cite{GK} uses a technique similar to ours and  
establishes the stability in the whole domain $D \times (-\delta, \delta)$
with more data. 
The main machinery 
is a Carleman estimate which requires information of the trace and the 
normal derivative of estimated function on the whole boundary 
$\ppp\OOO \times (-\delta, \delta)$.    
In this article, we apply a different Carleman estimate and prove 
conditional stability estimates in some subdomain of 
$D \times (-\delta,\delta)$.

We remark that we need not any data on $D \times \{\ell\} \times 
(-\delta, \delta)$.
Our approach is new for the inverse problem of determining functions which 
are independent of one component of the spatial variables.
Our formulation of the inverse problem is for example motivated by
the following.  Choosing the $x_n$-axis along the depth, we would like to
determine functions which are independent of the depth variable without
any data on the bottom $x_n=\ell$ of the cylindrical domain, 
but only data on the surface $x_n=0$ and
the side boundary data which can be approximated by
$\ppp_{x_n}u(\xi_k, x_n,t)$, $0<x_n<\ell$, $-\delta<t<\delta$ with
fixed probe points $\xi_1, ..., \xi_N \in \Gamma \subset \ppp D$.
Our main stability results guarantee that such data becomes more accurate
as the number $N$ of probe holes $\{ (\xi_k, x_n) \in \R^n;\, 0<x_n<\ell \}$,
$k=1,..., N$, increases. 

More precisely for the proof, we adapt the approach
by Bukhgeim and Klibanov \cite{BK} and Klibanov \cite{Kl2}
which originally was  to establish  the unique determination  and stability 
a $t$-independent source terms and coefficients of evolution equations.  In our case we determine an
$x_n$-independent unknown function.  As for inverse problems
by Carleman estimates, see Beilina and Klibanov \cite{BeiKl},
Bellassoued and Yamamoto \cite{BY}, Imanuvilov and Yamamoto
\cite{IY1998, IY1}, Klibanov and Timonov \cite{KT},
Yamamoto \cite{Ya} and the references therein.  Moreover we can refer
for example to Chapter 3, Section 3 in Isakov \cite{Is}
as for related inverse problems of determining
functions which are independent of one component of variables.

This paper is composed of three sections and one appendix.  
In Section 2, we establish a key
Carleman estimate and in Section 3, we complete the proof of Theorems 1.1 and 
1.2 and Corollary 1.1.

\section{Key Carleman estimate}

We set 
$$
\OOO_{\pm} := D \times (-\ell,\ell), \quad 
Q_{\pm} := \OOO_{\pm} \times (-\delta, \delta).
$$
We recall that $D \subset \R^{n-1}$ is a bounded domain with smooth boundary 
$\ppp D$, and $\OOO = D \times (0,\ell)$, 
$Q = \OOO  \times (-\delta,\delta)$.

We start the proof of Carleman estimate with construction of the weight 
function.
We can prove (e.g., Imanuvilov \cite{Ima1995}):
\\
{\bf Lemma 2.1.}\\{\it 
For a given subdomain $D_0 \subset D$ satisfying (1.5), there exists function
$d \in C^2(\ooo{D})$ such that
$$
\left\{ \begin{array}{rl}
& d(x') \ge 0 \quad \mbox{for $x' \in \ooo{D}$}, \\
& d(x') > 0 \quad \mbox{for $x' \in \ooo{D_0}$}, \\
& d(x') = 0 \quad \mbox{for $x' \in \ppp D \setminus \Gamma$}, \quad
\vert \nabla d(x') \vert > 0 \quad \mbox{for $x' \in \ooo{D}$}.
\end{array}\right.                   \eqno{(2.1)}
$$}
{\bf Proof of Lemma 2.1.}\\
For $\Gamma \subset \ppp D$, we choose a bounded domain $E$
with smooth boundary such that
$$
D \subsetneqq E, \quad \ooo{\Gamma} = \ooo{\ppp D\cap E},
\quad \ppp D\setminus\Gamma \subset \ppp E.
$$
In particular, $E\setminus\ooo{D}$ contains some non-empty
open subset.
We note that $E$ can be constructed as the interior of a union
of $\ooo{D}$ and the closure of a non-empty domain
$\widehat{D}$ satisfying $\widehat{D} \subset \overline{\R^3
\setminus D}$ and $\ppp\widehat{D} \cap \ppp D = \Gamma$.

We choose a domain $\omega$ such that
$\ooo{\omega} \subset E \setminus \ooo{D}$.
Then, by  \cite[Lemma 1.2]{Ima1995} (see also \cite[Lemma 2.1]{IY1998}), 
we can find $d\in C^2(\ooo{E})$ such that
$$
d>0 \quad \mbox{in $E$}, \quad |\nabla d| > 0 \quad \mbox{on
$\ooo{E\setminus \omega}$}, \quad d=0 \quad \mbox{on $\ppp E$}.
$$
This $d$ is our desired function, and the proof of Lemma 2.1 is complete. 
                                        $\blacksquare$
\\

We set
$$
\psi(x,t) = d(x') - \alpha x_n^2 - \beta t^2,
\quad \va(x,t) = e^{\la\psi(x,t)}
$$
where the positive constants $\alpha, \beta$ are chosen later and
$\la$ is a sufficiently large constant.
\\

{\bf Lemma 2.2 (Carleman estimate).}\\
{\it 
Let $p_0 \in L^{\infty}(D \times (-\delta, \delta))$ be given.
Then there exist constants $C>0$ and $s_0>0$ such that
\begin{align*}
& \int_{Q_{\pm}} \left( \frac{1}{s} \sum_{i,j=1}^n \vert \ppp_{x_i}\ppp_{x_j}u
\vert^2 + s\vert \nabla u\vert^2 + s^3\vert u\vert^2\right) \weight dxdt \\
\le &C\biggl(\int_{Q_{\pm}} \vert \ppp_tu - \Delta u - p_0u\vert^2 \weight dxdt
+ s^3\int_{\ppp\OOO_{\pm} \times (-\delta,\delta)}
(\vert \nabla_{x,t} u\vert^2 + \vert u\vert^2) \weight d\sigma dt\\
+& \frac 1s \Vert u e^{s\varphi}\Vert^2
_{L^2(-\delta,\delta;H^2(\partial\Omega_{\pm}))}\\
+& s^3\int_{\OOO_{\pm}} (\vert \nabla u(x,\delta)\vert^2
+ \vert u(x,\delta)\vert^2 + \vert \nabla u(x,-\delta)\vert^2
+ \vert u(x,-\delta)\vert^2) e^{2s\va(x,\delta)} dx\biggr)
\end{align*}
for all $s>s_0$ and all $u \in H^{2,1}(Q_{\pm})\cap H^1(-\delta,\delta;
H^1(\OOO_{\pm}))$ satisfying $u\in L^2(-\delta,\delta;H^2(\ppp\OOOPM))$.   
Here $s_0$ can be uniformly
chosen if $\Vert p_0\Vert_{L^{\infty}(D \times (-\delta, \delta))}$ is
bounded.}
\\

The proof of Lemma 2.2 is based on a classical Carleman estimate and is given 
in Appendix for completeness.

\section{Proof of Main Results}

\subsection{Proof of Theorem 1.1 (i)}
We recall
$$
u(\xxz,t) = \ppp_{x_n}u(\xxz,t) = 0, \quad x' \in D, \, -\delta<t<\delta
                                                   \eqno{(3.1)}
$$
and
$$
\ppp_tu = \Delta u + p_0(x',t) u + R(\xxx,t)f(x',t), \quad x'\in D, \,
0 < x < \ell, \, -\delta<t<\delta.               \eqno{(3.2)}
$$
We make the even extension of functions $u$ and $R$ in the variable $x_n$:
$$
u(\xxx,t) =
\left\{\begin{array}{rl}
& u(\xxx,t), \quad x_n \ge 0, \\
& u(x',-x_n,t), \quad x_n < 0,
\end{array}\right.
\quad
R(\xxx,t) =
\left\{\begin{array}{rl}
&R
(\xxx,t), \quad x_n \ge 0, \\
&R
(x',-x_n,t), \quad x_n < 0.
\end{array}\right.
$$
By (3.1) we can verify
$$
\ppp_{x_n}u(\xxx,t) =
\left\{\begin{array}{rl}
& \ppp_{x_n}u(\xxx,t), \quad x_n \ge 0, \\
& -\ppp_{x_n}u(x',-x_n,t), \quad x_n < 0,
\end{array}\right.
$$
and
$$
\ppp_{x_n}^2 u(\xxx,t) =
\left\{\begin{array}{rl}
& \ppp_{x_n}^2u(\xxx,t), \quad x_n \ge 0, \\
& \ppp_{x_n}^2u(x',-x_n,t), \quad x_n < 0.
\end{array}\right.
$$
Hence we can prove that
$\ppp_{x_n}^3u \in L^2(D \times (-\ell,\ell) \times (-\delta, \delta))$ and so
$\ppp_{x_n}u \in L^2(-\delta, \delta; H^2(D \times (-\ell,\ell)))$, and
$\ppp_{x_n}u \in H^1(-\delta, \delta; H^1(D \times (-\ell,\ell))$.
Moreover
$$
\ppp_{x_n}R
(\xxx,t) =
\left\{\begin{array}{rl}
& \ppp_{x_n}R
(\xxx,t), \quad x_n \ge 0, \\
& -\ppp_{x_n}R(
x',-x_n,t), \quad x_n < 0,
\end{array}\right.
$$
and so $\ppp_{x_n}R \in L^{\infty}(D \times (-\ell,\ell)
\times (-\delta, \delta))$ by (1.3).

Thus (3.2) yields
$$\left\{ \begin{array}{rl}
& \ppp_tu = \Delta u + p_0u + R(x',x_n,t)f(x',t), \quad x'\in D,
\, -\ell<x_n<\ell,\, -\delta < t < \delta, \\
& u(\xxz,t) = \ppp_{x_n} u(\xxz,t) = 0, \quad x'\in D,\, -\delta<t<\delta.
\end{array}\right.
                                    \eqno{(3.3)}
$$
Set
$$
y = \ppp_{x_n}u.
$$
Then $y$ satisfies
$$\left\{ \begin{array}{rl}
& \ppp_ty = \Delta y + p_0y + \ppp_{x_n}R(x',x_n,t)f(x',t), \quad x'\in D, \,
-\ell<x_n<\ell,\, -\delta < t < \delta, \\
& y(\xxz,t) = 0, \quad x'\in D,\, -\delta<t<\delta.
\end{array}\right.                                   \eqno{(3.4)}
$$
\\

Now we will specify the weight function $\psi(x,t) : = d(x') - \alpha
x_n^2
- \beta t^2$ for the Carleman estimate Lemma 2.2.
We set
$$
d_0:= \min_{x'\in\ooo{D_0}} d(x'), \quad
d_1:= \max_{x'\in\ooo{D}} d(x').
$$
By Lemma 2.1, we see $d_0 > 0$.
We choose $\delta_0>0$ such that 
$$
\delta_0 < \left( \frac{d_0}{d_1}\right)^{\frac{1}{2}}\delta.  \eqno{(3.5)}
$$
We note that $0 < \delta_0 < \delta$.  Then, since (3.5) yields
$$
\frac{d_1-d_0}{\delta^2-\delta_0^2} < \frac{d_0}{\delta_0^2},
$$
we can choose $\beta > 0$ such that
$$
\frac{d_1-d_0}{\delta^2-\delta_0^2} < \beta < \frac{d_0}{\delta_0^2}.
                                                   \eqno{(3.6)}
$$
Finally choose $\alpha > 0$ sufficiently large such that
$$
d_1 - d_0 + \beta \delta_0^2 < \alpha \ell^2.           \eqno{(3.7)}
$$
Then inequalities (3.5) - (3.7) imply 
$$
\left\{ \begin{array}{rl}
& d_1 - \beta\delta^2 < d_0 - \beta \delta_0^2, \\
& 0 < d_0 - \beta \delta_0^2, \\
& d_1 - \alpha\ell^2 < d_0 - \beta \delta_0^2.
\end{array}\right.                                 \eqno{(3.8)}
$$
Here we recall 
$$
\psi(x,t) = d(x') - \alpha x_n^2 - \beta t^2,
\quad \va(x,t) = e^{\la\psi(x,t)}.
$$

Inequalities (3.8) imply
\begin{align*}
& \max\biggl\{
\max_{x\in \ooo{\OOO}} \psi(x,\delta), \,
\max_{x'\in \ppp D \setminus \Gamma, \, -\ell\le x_n \le \ell, \,
-\delta\le t\le \delta} \psi(x,t),
\max_{x'\in \ooo{D},\, -\delta\le t\le \delta} \psi(x,\pm\ell,t)
\biggr\} \\
< & \min_{x'\in \ooo{D_0},\, -\delta_0\le t \le \delta_0} \psi(\xxz,t)
= d_0 - \beta \delta_0^2.
\end{align*}
Therefore
$$
\sigma_1 := \max\biggl\{
\max_{x\in \ooo{\OOO}} \va(x,\delta), \,
\max_{x'\in \ppp D \setminus \Gamma, \, -\ell\le x_n \le \ell, \,
-\delta\le t\le \delta} \va(x,t),
\max_{x'\in \ooo{D},\, -\delta\le t\le \delta} \va(x,\pm\ell,t)
\biggr\}
$$
$$
< \sigma_0 :=
\min_{x'\in \ooo{D_0},\, -\delta_0\le t \le \delta_0} \va(\xxz,t).
                                            \eqno{(3.9)}
$$

Next in terms of (3.9), we estimate the integral
\begin{align*}
&\int_{D\times (-\delta,\delta)} \vert \ppp_{x_n}y(\xxz,t)\vert^2
e^{2s\va(\xxz,t)} dx' dt\\
= &\int^0_{-\ell} \frac{d}{d\xi}
\left( \int_{D\times (-\delta,\delta)}\vert \ppp_{x_n}y(x',\xi,t)\vert^2
e^{2s\va(x',\xi,t)} dx' dt\right) d\xi\\
+ & \int_{D\times (-\delta,\delta)} \vert \ppp_{x_n}y(x',-\ell,t)\vert^2
e^{2s\va(x',-\ell,t)} dx'dt \\
= & \int^0_{-\ell} \int_{D\times (-\delta,\delta)}
(\vert \ppp_{x_n}y(x',\xi,t)\vert^2 2s(\ppp_{\xi}\va)(x',\xi,t)\\
+ & 2(\ppp_{x_n}y)(x',\xi,t)(\ppp_{x_n}^2y)(x',\xi,t) 
e^{2s\va(x',\xi,t)} dx' dt d\xi
                                                               \\
+ & \int_{D\times (-\delta,\delta)} \vert \ppp_{x_n}y(x',-\ell,t)\vert^2
e^{2s\va(x',-\ell,t)} dx'dt.
\end{align*}
By (3.9) and (1.6), we have
$$
\max_{x'\in\ooo{D},-\delta\le t \le \delta} \va(x',-\ell,t) \le \sigma_1
$$
and $\Vert \ppp_{x_n}y(\cdot,-\ell,\cdot)\Vert_{L^2(D \times (-\delta,\delta))}
\le M$.
Therefore
$$
\int_{D\times (-\delta,\delta)} \vert \ppp_{x_n}y(x',0, t)\vert^2
e^{2s\va(x',0, t)} dx' dt                               \eqno{(3.10)}
$$
$$
\le C\int_{Q_{\pm}} (s\vert \ppp_{x_n}y\vert^2 + \vert \ppp_{x_n}y\vert
\vert \ppp_{x_n}^2y\vert)(x,t) \weight dxdt
+ CM^2e^{2s\sigma_1}.
$$
Since 
$$
\vert \ppp_{x_n}^2y\vert \vert \ppp_{x_n}y\vert 
\le \frac{1}{2} \left( \frac{1}{s}\vert \ppp_{x_n}^2y\vert^2 
+ s\vert \ppp_{x_n}y\vert^2\right)
$$
and 
$\partial_{x_n}^2u = \ppp_{x_n}y$, we rewrite the first term in estimate
(3.10) to have
$$
\int_{D\times (-\delta,\delta)} \vert \ppp_{x_n}^2u(x', 0,t)\vert^2
e^{2s\va(x',0,t)} dx'dt
$$
$$
\le C\int_{Q_{\pm}} 
\left(s\vert \ppp_{x_n}y\vert^2 + \frac{1}{s}\vert \ppp_{x_n}^2y
\vert^2 \right) \weight dxdt + CM^2e^{2s\sigma_1}        \eqno{(3.11)}
$$
We apply Lemma 2.2 to system (3.4), and we obtain
$$
\int_{Q_{\pm}} \left(\frac {1}{s} \vert \ppp_{x_n}^2y\vert^2
+ s\vert \ppp_{x_n}y\vert^2 \right) \weight dxdt
                                                            \eqno{(3.12)}
$$
\begin{align*}
\le & C\int_{Q_{\pm}} \vert f(x',t)\vert^2 e^{2s\va(\xxx,t)} dx'dx_ndt\\
+ & Cs^3 \int_{\ppp\OOO_{\pm} \times (-\delta,\delta)}
(\vert \nabla_{x,t}y\vert^2 + \vert y\vert^2)
\weight d\sigma dt 
+ C\Vert y e^{s\varphi}\Vert^2_{L^2(-\delta,\delta;H^2(\partial\OOOPM))}\\
+ & Cs^3\int_{\OOOPM} (\vert \nabla y (x,\delta)\vert^2
+ \vert y(x,\delta)\vert^2 + \vert \nabla y (x,-\delta)\vert^2
+ \vert y(x,-\delta)\vert^2) ee^{2s\va(x,\delta)} dx.
\end{align*}
Here we have
$$
\int_{\ppp\OOO_{\pm} \times (-\delta,\delta)}
(\vert \nabla_{x,t}y\vert^2 + \vert y\vert^2)
\weight d\sigma dt                        \eqno{(3.13)}
$$
\begin{align*}
=& \biggl( \int_{\Gamma \times (-\ell,\ell) \times (-\delta, \delta)}
+ \int_{(\ppp D \setminus \Gamma) \times (-\ell,\ell)
\times (-\delta,\delta)}
+ \int_{D\times \{\ell\} \times (-\delta,\delta)} \\
+ & \int_{D\times \{-\ell\} \times (-\delta,\delta)} \biggr)
(\vert \nabla_{x,t}y\vert^2 + \vert y\vert^2)
\weight d\sigma dt 
\le Ce^{Cs}\mathcal{D}(u)^2 + Ce^{2s\sigma_1}M^2
\end{align*}
by (1.6), (1.7) and (3.9).  Moreover,
\begin{align*}
& \Vert ye^{s\va}\Vert^2_{L^2(-\delta,\delta;H^2(\ppp\OOO_{\pm}))}
= \Vert ye^{s\va}\Vert^2_{L^2(-\delta,\delta;H^2(\Gamma \times (-\ell,\ell))}
+ \Vert ye^{s\va}\Vert^2_{L^2(-\delta,\delta;H^2((\ppp D \setminus \Gamma)
\times (-\ell,\ell))}\\
+& \Vert ye^{s\va}\Vert^2_{L^2(-\delta,\delta;H^2(D \times \{\ell\}))}
+ \Vert ye^{s\va}\Vert^2_{L^2(-\delta,\delta;H^2(D \times \{-\ell\}))}.
\end{align*}
We can directly verify that 
$$
\Vert ye^{s\va}(\cdot,t)\Vert^2_{H^2(\gamma)}
\le Cs^4\int_{\gamma} \left( \sum_{i,j=1}^n \vert \ppp_{x_i}\ppp_{x_j}y\vert^2
+ \vert \nabla y\vert^2 + \vert y\vert^2\right) e^{2s\va(x,t)} d\sigma,
$$
where $\gamma = \Gamma \times (-\ell, \ell)$ or 
$= (\ppp D \setminus \Gamma) \times (-\ell, \ell)$ or 
$= D \times \{ \ell\}$ or $= D \times \{ -\ell\}$.  Therefore, again using
(1.6), (1.7) and (3.9), we obtain
$$
\Vert ye^{s\va}\Vert^2_{L^2(-\delta,\delta;H^2(\ppp\OOO_{\pm}))}
\le Cs^4e^{Cs}\mathcal{D}(u)^2 + Cs^4e^{2s\sigma_1}M^2    \eqno{(3.14)}
$$
for all large $s>0$.
Since 
\begin{align*}
& \int_{\OOOPM} (\vert \nabla y(x,\delta)\vert^2 
+ \vert y(x,\delta)\vert^2 + \vert \nabla y(x,-\delta)\vert^2
+ \vert y(x,-\delta) \vert^2) e^{2s\va(x,\delta)} 
dx                              \\
\le& CM^2e^{2s\sigma_1}
\end{align*}
by (1.6) and (3.9), we substitute (3.13) and (3.14) into (3.12), and reach
$$
\int_{Q_{\pm}} \left(\frac{1}{s}\vert \ppp_{x_n}^2y\vert^2 
+ s\vert \ppp_{x_n}y \vert^2 \right) \weight dxdt 
                                                     \eqno{(3.15)}
$$
$$
\le C\int_{Q_{\pm}} \vert f(x',t)\vert^2 e^{2s\va(\xxx,t)} dx'dx_ndt
+ Cs^4e^{Cs}\mathcal{D}(u)^2 + Cs^4e^{2s\sigma_1}M^2
$$
for all large $s>0$.

From (3.11) and (3.15), we have
\begin{align*}
&\int_{D\times (-\delta,\delta)} \vert \ppp_{x_n}^2u(x', 0,t)\vert^2
e^{2s\va(x',0,t)} dx'dt                        \\
\le & C\int_{D\times (-\delta,\delta)}\vert f(x',t)\vert^2 e^{2s\va(\xxz,t)}
\left( \int^{\ell}_{-\ell} e^{2s(\va(\xxx,t) - \va(\xxz,t))} dx_n\right)
                                                         dx'dt
\end{align*}
$$
+ Cs^4e^{Cs}\mathcal{D}(u)^2 + Cs^4e^{2s\sigma_1}M^2.        \eqno{(3.16)}
$$
Since $f(x',t)R(\xxz,t) = - \ppp_{x_n}^2u(\xxz,t)$ for $x' \in D$ and
$-\delta<t<\delta$ by (3.3), from (3.16) and (1.4), we obtain
\begin{align*}
&\int_{D\times (-\delta,\delta)} \vert f(x',t)\vert^2
e^{2s\va(\xxz,t)} dx' dt
\le C\int_{D\times (-\delta,\delta)} \vert \ppp_{x_n}^2u(x',0,t)\vert^2
e^{2s\va(\xxz,t)} dx' dt\\
\le &C\int_{D\times (-\delta,\delta)} \vert f(x',t)\vert^2
e^{2s\va(\xxz,t)} dx' dt
\left( \int^{\ell}_{-\ell} e^{2s(\va(\xxx,t) - \va(\xxz,t))} dx_n\right)
dx'dt\\
+ &Cs^4e^{Cs}\mathcal{D}(u)^2 + Cs^4e^{2s\sigma_1}M^2.
\end{align*}
Here we have
\begin{align*}
& e^{2s(\va(\xxx,t) - \va(\xxz,t))}
= e^{2se^{\la(d(x')-\beta t^2)}(e^{-\la\alpha x_n^2}-1)}\\
= & e^{-2se^{\la(d(x')-\beta t^2)}(1-e^{-\la\alpha x_n^2})}
\le e^{-2sc_0(1-e^{-\la\alpha x_n^2})},
\end{align*}
where
$$
c_0 := \min_{x'\in \ooo{D}, -\delta\le t \le \delta}
e^{\la(d(x') - \beta t^2)}.
$$
Therefore the Lebesgue theorem yields
$$
\int^{\ell}_{-\ell} e^{2s(\va(\xxx,t) - \va(\xxz,t))} dx_n
\le \int^{\ell}_{-\ell} e^{-2sc_0(1-e^{-\la\alpha x_n^2})}dx_n
= o(1)\quad \mbox{as}\,\, s \to \infty.
$$
Hence
\begin{align*}
& \int_{D\times (-\delta,\delta)} \vert f(x',t)\vert^2
e^{2s\va(\xxz,t)} dx' dt\\
= &o(1)\int_{D\times (-\delta,\delta)} \vert f(x',t)\vert^2
e^{2s\va(\xxz,t)} dx' dt
+ Cs^4e^{Cs}\mathcal{D}(u)^2 + Cs^4e^{2s\sigma_1}M^2.
\end{align*}
Choosing the parameter $s>0$ large, 
we can absorb the first term on the right-hand side
into the left-hand side, and
$$
\int_{D\times (-\delta,\delta)} \vert f(x',t)\vert^2
e^{2s\va(\xxz,t)} dx' dt
\le Cs^4e^{Cs}\mathcal{D}(u)^2 + Cs^4e^{2s\sigma_1}M^2
$$
for all sufficiently large $s>0$.  Shrinking $D \times (-\delta,\delta)$ to
$D_0 \times (-\delta_0,\delta_0)$ in the integral on the left-hand side
and applying the definition of $\sigma_0$ in (3.9), we obtain
$$
\Vert f\Vert^2_{L^2(D_0\times (-\delta_0,\delta_0))}e^{2s\sigma_0}
\le Cs^4e^{Cs}\mathcal{D}(u)^2 + Cs^4e^{2s\sigma_1}M^2,
$$
that is,
$$
\Vert f\Vert^2_{L^2(D_0\times (-\delta_0,\delta_0))}
\le Cs^4e^{Cs}\mathcal{D}(u)^2 + Cs^4e^{-2s(\sigma_0-\sigma_1)}M^2
$$
for all $s>s_0$: some constant.  
By $\sigma_0-\sigma_1 > 0$, we see that $\sup_{s\ge 0}
s^4e^{-s(\sigma_0-\sigma_1)} < \infty$, and choosing  
a constant $C_1>0$ satisfying $s^4e^{Cs} \le e^{C_1s}$ for
all $s\ge 0$, we have
$$
\Vert f\Vert^2_{L^2(D_0\times (-\delta_0,\delta_0))}
\le Ce^{C_1s}\mathcal{D}(u)^2 + C_2e^{-s(\sigma_0-\sigma_1)}M^2
$$
for all $s\ge s_0$.  Replacing $C$ by $Ce^{C_1s_0}$ and changing  
$s$ into $s+s_0$ with $s\ge 0$, we obtain
$$
\Vert f\Vert^2_{L^2(D_0\times (-\delta_0,\delta_0))}
\le Ce^{C_1s}\mathcal{D}(u)^2 + C_2e^{-s(\sigma_0-\sigma_1)}M^2
                                             \eqno{(3.17)}
$$
for all $s\ge 0$.
We choose $s>0$ in order that the right-hand side of (3.17) is small and 
consider the two cases separately.
\\
{\bf Case 1: $M > \mathcal{D}(u)$}:\\
We choose $s>0$ such that 
$$
e^{C_1s}\mathcal{D}(u)^2 = e^{-s(\sigma_0-\sigma_1)}M^2, 
\quad \mbox{that is,} \quad 
s = \frac{2}{C_1+\sigma_0-\sigma_1}\log \frac{M}{\mathcal{D}(u)} > 0.
$$
Substituting this value of $s$ into (3.17), we reach
$$
\Vert f\Vert^2_{L^2(D_0\times (-\delta_0,\delta_0))}
\le 2M^{2(1-\theta)}\mathcal{D}(u)^{2\theta},
$$
where $\theta = \frac{\sigma_0-\sigma_1}{C_1+\sigma_0-\sigma_1} 
\in (0,1)$.
\\
{\bf Case 2: $M \le \mathcal{D}(u)$}:\\
Setting $s=0$ in (3.17), we directly obtain
$$
\Vert f\Vert^2_{L^2(D_0\times (-\delta_0,\delta_0))}
\le (C+C_2)\mathcal{D}(u)^2.
$$
By the definitions (1.6) and (1.7) of $M$ and $\mathcal{D}(u)$, we see
$\mathcal{D}(u) \le CM$, and so 
$$
\Vert f\Vert^2_{L^2(D_0\times (-\delta_0,\delta_0))}
\le (C+C_2)\mathcal{D}(u)^{2\theta}M^{2(1-\theta)}.
$$
Thus the proof of Theorem 1.1 (i) is complete. $\blacksquare$

\subsection{Proof of Theorem 1.1 (ii)}

We fix  an arbitrary  point $x_0'$ from an interior of the set  $\Gamma.$ 
For sufficiently small $\ep > 0$, we choose
$$
\www{D} = D \cap \{ x' \in \R^{n-1};\, \vert x' - x_0'\vert < 2\ep\},
\quad
D_1 = D \cap \{ x' \in \R^{n-1};\, \vert x' - x_0'\vert < \ep\}
$$
such that $\www{D} \cap \ppp D \subsetneq \Gamma$.
Then, for small $\ep > 0$, replacing $D_0$ and $D$ respectively
by $D_1$ and $\www{D}$, in terms of Lemma 2.1 we can construct 
$d \in C^2(\ooo{\www{D}})$ satisfying (2.1) with $\www{D}$ replacing 
$D$.  Let $\delta_1 > 0$ be chosen arbitrarily such that
$0 < \delta_1 < \delta$.  Then for sufficiently small $\ep > 0$, we can
verify
$$
\left( \frac{\delta_1}{\delta}\right)^2
< \frac{ \min_{x'\in \ooo{D_1}} d(x')}
{\max_{x'\in \ooo{\www{D}}} d(x')} < 1.             \eqno{(3.18)}
$$
This is possible because
$$
\lim_{\ep \to 0} \frac{ \min_{x'\in \ooo{D_1}} d(x')}
{\max_{x'\in \ooo{\www{D}}} d(x')} = 1
$$
and $\frac{\delta_1}{\delta} < 1$.
Replacing $D_0$ and $D$ by $D_1$ and $\www{D}$ respectively, we apply
(3.18) instead of (3.5), and argue similarly to the proof of Theorem 1.1 (i),
so that the proof of Theorem 1.1 (ii) is complete. $\blacksquare$

\subsection{Proof of Corollary 1.1}

By Theorem 1.1 (i) and its proof in Section 3.1, for arbitrary
subdomain $D_0 \subset D$ satisfying (1.4), we see that
$f=0$ in $D_0 \times (-\delta_0,\delta_0)$ where $\delta_0>0$ satisfies
(1.8).  Therefore the trace theorem yields
$f=0$ in $D_0 \times \{0\}$.  Since $D_0 \subset D$ can be chosen
arbitrarily provided that (1.5) holds, we see that
$f=0$ in $D \times \{0\}$.  Thus the proof of the corollary is complete.
$\blacksquare$

{\bf Proof of Theorem 1.2.}
The theorem follows directly from Theorem 1.1.  Indeed, setting 
$u := v(p) - v(q)$,
$p_0 := p$, $R:= v(q)$ and $f:= p-q$ and taking the difference
between the two corresponding equations with $v(p)$ and $v(q)$, 
we reduce Theorem 1.2 to Theorem 1.1.
                                             $\blacksquare$

\section*{Appendix. Proof of Lemma 2.2}

The following inequality is a standard Carleman estimate:
There exist constants $C>0$ and $s_0>0$ such that
$$
\int_{Q_{\pm}} \left( \frac{1}{s}\vert \Delta u\vert^2 
+ s\vert \nabla u\vert^2 + s^3\vert u\vert^2 \right) \weight dxdt 
                                        \eqno{(1)}
$$
\begin{align*}
\le &C\biggl(\int_{Q_{\pm}} \vert \ppp_tu - \Delta u - p_0u\vert^2 \weight dxdt
+ s^3\int_{\ppp\OOO_{\pm} \times (-\delta,\delta)}
(\vert \nabla_{x,t} u\vert^2 + \vert u\vert^2) \weight d\sigma dt\\
+ & s^3\int_{\OOO_{\pm}} (\vert \nabla u(x,\delta)\vert^2
+ \vert u(x,\delta)\vert^2 + \vert \nabla u(x,-\delta)\vert^2
+ \vert u(x,-\delta)\vert^2) e^{2s\va(x,\delta)} dx\biggr)
=:J
\end{align*}
for all $s>s_0$.  See e.g., \cite{BY}, \cite{Ya} as for the proof.

We recall that $\OOO_{\pm} = D \times (-\ell,\ell)$.
Henceforth $\nu = (\nu_1, ..., \nu_n)$ denotes the unit outward normal vector 
to $\ppp\OOOPM$.
Then we prove

{\bf Lemma 1.}\\
{\it 
Let $w \in H^2(\OOO_{\pm})$ satisfy $w\in H^2(\ppp\OOO_{\pm})$.
Then
$$
\sum_{i,j=1}^n \Vert \ppp_{x_i}\ppp_{x_j}w\Vert^2_{L^2(\OOOPM)} 
+ \sum_{i,j=1}^n \int_{\ppp\OOOPM} \ppp_{x_i}w((\ppp_{x_j}\ppp_{x_j}w)\nu_i
- (\ppp_{x_i}\ppp_{x_j}w)\nu_j) d\sigma = \Vert \Delta w\Vert^2_{L^2(\OOOPM)}.
$$
In particular,
$$
\sum_{i,j=1}^n \Vert \ppp_{x_i}\ppp_{x_j}w\Vert_{L^2(\OOOPM)}
\le C(\Vert \Delta w\Vert_{L^2(\OOOPM)} + \Vert w\Vert_{H^2(\ppp\OOOPM)}),
$$
where the constant $C>0$ depends only on $\OOOPM$.
}
\\
{\bf Remark.}  We can relax the norm $\Vert w\Vert_{H^2(\ppp\OOOPM)}$ by
$\Vert w\Vert_{H^{\frac{3}{2}}(\ppp\OOOPM)}$, but this choice of the norm
of Dirichlet data is sufficient for the proof of the main theorems.
\\
{\bf Proof.}\\
By the density argument, it suffices to assume that $w \in 
C^{\infty}(\ooo{\OOOPM})$.
Then the integration by parts yield
\begin{align*}
& \int_{\OOOPM} \vert \Delta w\vert^2 dx
= \sum_{i,j=1}^n \int_{\OOOPM} (\ppp_{x_i}\ppp_{x_i}w)
(\ppp_{x_j}\ppp_{x_j}w) dx \\
=& -\sum_{i,j=1}^n \int_{\OOOPM} (\ppp_{x_i}w)(\ppp_{x_j}\ppp_{x_i}
\ppp_{x_j}w) dx
+ \sum_{i,j=1}^n \int_{\ppp\OOOPM} (\ppp_{x_i}w)(\ppp_{x_j}\ppp_{x_j}w)
\nu_i d\sigma\\
=& \sum_{i,j=1}^n \int_{\OOOPM} \vert \ppp_{x_i}\ppp_{x_j}w\vert^2 dx
+ \sum_{i,j=1}^n \int_{\ppp\OOOPM} 
\ppp_{x_i}w((\ppp_{x_j}\ppp_{x_j}w)\nu_i - (\ppp_{x_i}\ppp_{x_j}w)\nu_j) 
d\sigma.
\end{align*}
Thus the proof of Lemma 1 is complete.    $\blacksquare$
\\

We arbitrarily fix $t \in (-\delta,\delta)$.
Direct calculations yield 
$$
\Delta (u(x,t)e^{s\va(x,t)})
= (\Delta u)e^{s\va} + 2s(\nabla u\cdot \nabla\va) e^{s\va}
+ (s^2\vert \nabla \va\vert^2 + s\Delta\va)ue^{s\va}, \quad 
x \in \OOOPM.
$$
Applying Lemma 1 to $u(\cdot,t)e^{s\va(\cdot,t)}$ and integrating over $t \in
(-\delta, \delta)$, we have
\begin{align*}
& \frac{1}{s} \sum_{i,j=1}^n \Vert \ppp_{x_i}\ppp_{x_j}
(ue^{s\va})\Vert^2_{L^2(Q_{\pm})} \\
\le& C\left( \frac{1}{s}\Vert (\Delta u)e^{s\va}\Vert^2_{L^2(Q_{\pm})}
+ s\Vert (\nabla u) e^{s\va}\Vert^2_{L^2(Q_{\pm})}
+ s^3\Vert u\Vert^2_{L^2(Q_{\pm})}\right)
+ \frac{C}{s}\Vert ue^{s\va}\Vert^2_{L^2(-\delta,\delta;H^2(\ppp\OOOPM))}.                        
\end{align*}
Applying (1), we obtain
$$
\frac{1}{s} \sum_{i,j=1}^n \Vert \ppp_{x_i}\ppp_{x_j}
(ue^{s\va})\Vert^2_{L^2(Q_{\pm})} 
\le CJ + \frac{C}{s}\Vert ue^{s\va}\Vert^2
_{L^2(-\delta,\delta;H^2(\ppp\OOOPM))}.
                                     \eqno{(2)}
$$                        
Since 
\begin{align*}
&\ppp_{x_i}\ppp_{x_j}(ue^{s\va})
= (\ppp_{x_i}\ppp_{x_j}u)e^{s\va}\\
+& se^{s\va}((\ppp_{x_i}u)\ppp_{x_j}\va + (\ppp_{x_j}u)\ppp_{x_i}\va)
+ (s^2 (\ppp_{x_i}\va)\ppp_{x_j}\va + s(\ppp_{x_i}\ppp_{x_j}\va))
e^{s\va}u,
\end{align*}
again using (1) and (2), we have
\begin{align*}
& \frac{1}{s} \sum_{i,j=1}^n \Vert (\ppp_{x_i}\ppp_{x_j}u)e^{s\va}
\Vert^2_{L^2(Q_{\pm})}\\
\le & C\left( \frac{1}{s} \sum_{i,j=1}^n \Vert \ppp_{x_i}\ppp_{x_j}
(ue^{s\va})\Vert^2_{L^2(Q_{\pm})}
+ s\Vert (\nabla u) e^{s\va}\Vert^2_{L^2(Q_{\pm})}
+ s^3\Vert u\Vert^2_{L^2(Q_{\pm})}\right)\\
\le & CJ + \frac{C}{s}\Vert ue^{s\va}\Vert^2
_{L^2(-\delta,\delta;H^2(\ppp\OOOPM))}. 
\end{align*}
Integrating over $t \in (-\delta, \delta)$, we complete the proof
of Lemma 2.2.                       $\blacksquare$

\section*{Acknowledgments}
The first author was partially supported by NSF grant DMS 1312900 and
Grant-in-Aid for Scientific Research (S)
15H05740 of Japan Society for the Promotion of Science.
The second author is partially supported by the French National Research 
Agency ANR (project MultiOnde) grant ANR-17-CE40-0029.
The third author was supported by Grant-in-Aid for Scientific Research (S)
15H05740 of Japan Society for the Promotion of Science and
by The National Natural Science Foundation of China
(no. 11771270, 91730303).
This work was prepared with the support of the "RUDN University Program 5-100".


\end{document}